\documentclass[12pt]{article}
\usepackage{amsmath,amsfonts}
\begin{document}

\baselineskip 20pt
\date{}
\title{On a boundary property of Blaschke products}

\author{Arthur A.~Danielyan and Spyros Pasias}

\maketitle 

\begin{abstract}
 
\noindent A Blaschke product has no radial limits on a subset $E$
of the unit circle $T$ but has unrestricted limit at each point of $T \setminus E$ if and only if $E$ is a closed set of
measure zero.


\end{abstract}

\begin{section}{Introduction.}

Below $\Delta$ and $T$ denote the open unit disk and the unit circle in the complex plane,
respectively. Let $m$ be the Lebesgue measure on $T$. 
The Blaschke products are forming
an important subclass of the well-known space $H^\infty$ of all bounded analytic functions in $\Delta$.
In this paper we consider a question on the boundary behavior of Blashke products.

Let $f$ be a function defined in $\Delta$ and let $z_0\in T$. Except the radial and non-tangential 
(angular) limits of $f$ at $z_0$, we also consider the unrestricted limit  of $f$ at $z_0$.
Recall that the unrestricted limit of $f$ at $z_0$ is the limit of $f$ when $z\in \Delta$ approaches to $z_0$
arbitrarily, if the limit exists. 

By a known theorem of Lindel{\"o}f, if
$f \in H^\infty$ then the existence of the radial limit at $z_0$ implies the existence of the non-tangential limit at the same point. 
But, of course, the unrestricted limit of $f$ at $z_0$ does not have to exist if just the radial limit exists at $z_0$.

The following theorem of Fatou is fundamental in the theory of boundary behavior of analytic functions.

\vspace{0.25 cm}

{\bf Theorem A.}\ {\it Let $f \in H^\infty$. Then the radial limits of $f$ exist  on $T$ except perhaps for a subset $E$ of measure zero.}
 
 \vspace{0.25 cm} 
 
 For Blaschke products, F. Riesz has proved the following more precise (than Theorem A) result. 
  
 \vspace{0.25 cm}

{\bf Theorem B.}\ {\it The moduli of the radial limits of a Blashke product $B$ are equal to 1 a.e. on $T$.}
 
 \vspace{0.25 cm}

 It is well known that the exceptional set $E$ in Theorem A is a $G_{\delta \sigma}$ set. But if we assume that 
 the function $f$ has unrestricted limits at the points of the set $T \setminus E$, then $E$ necessarily becomes 
 an $F_\sigma $ set (in fact this statement is obvious; cf \cite{Dan1}). 
 The converse statement was proved in \cite{Dan1} on the base of the method of \cite{Dan}, and the 
 result can be formulated as follows.




\vspace{0.25 cm}

{\bf Theorem C.}\ {\it  Let $E$ be a set on $T$. Then there exists a function $f \in H^\infty$ which has no radial limits on $E$
but has unrestricted limit at each point of $T \setminus E$ if and only if $E$ is an $F_\sigma$  set of measure zero.}
 
\vspace{0.25 cm}

Note that the proof of Theorem C is elementary and 
its sufficiency part implies
a well-known
theorem of Lohwater and Piranian \cite{LoPi} 
as an obvious corollary.



The purpose of the present paper is proving the analogues of Theorem C for Blaschke products.  
The result is the following:

 \vspace{0.25 cm}

{\bf Theorem 1.}\ {\it  Let $E$ be a set on $T$. Then there exists a Blaschke product which has no radial limits on $E$
but has unrestricted limit at each point of $T \setminus E$ if and only if $E$ is a closed set of measure zero.}
 
\vspace{0.25 cm}


The proof of Theorem 1 uses in particular Theorem C and some results 
on the boundary behavior of Blaschke products 
due to R. D. Berman \cite{Ber} and A. Nicolau \cite{Nic}. 

\end{section}

\begin{section}{Some auxiliary results.}  

We follow the presentation of Nicolau \cite{Nic} to formulate the results of
 Berman \cite{Ber} and Nicolau \cite{Nic}, respectively. The next theorem is due to Berman (cf. \cite{Nic}, p. 250).

  \vspace{0.25 cm}

{\bf Theorem D.}\ {\it Let E be a subset of the unit circle of zero Lebesgue
measure and of type $F_\sigma$ and $G_\delta$. Then there exist Blaschke products $B_0$ and $B_1$
such that:

(i) $B_0$ extends analytically to $T \setminus \overline E $ and $\lim_{r \rightarrow 1}B_0(re^{it})=0$ if and only if $e^{it} \in E$;

(ii)  $\lim_{r \rightarrow 1}B_1(re^{it})=1$ if and only if $e^{it} \in E$.}
 
\vspace{0.25 cm}

Now we formulate a theorem of Nicolau (see \cite{Nic}, Proposition on p. 251).
 
  \vspace{0.25 cm}

{\bf Theorem E.}\ {\it Let E be a subset of the unit circle. Assume that there exist a Blaschke products $B_0$ 
that extends analytically to $T\setminus \overline E$ with $\lim_{r \rightarrow 1}B_0(re^{it})=0$ for $e^{it} \in E$,
and an analytic function $f_1$ in the unit ball of $H^\infty$, $f_1 \not \equiv 1,$ 
such that $\lim_{r \rightarrow 1}f_1(re^{it})=1$ for $e^{it} \in E$. 
Then for each analytic function $g$ in the unit ball of $H^\infty$, there exists a Blaschke product $I$ 
 that extends analytically to $T \setminus \overline E$, such that
 
 $$\lim_{r \rightarrow 1} [I(re^{it}) -  g(re^{it})]=0 \  \ {\textit for}  \  \    e^{it} \in E.$$}
 

In this result, of course, the set $E$ in fact is of Lebesgue measure zero since if $m(E)>0$, then already 
the condition
 $\lim_{r \rightarrow 1}B_0(re^{it})=0$ for $e^{it} \in E$ will imply that $B_0 \equiv 0$, which is impossible as $B_0$ is a Blaschke product.

As noted in \cite{Nic}, p. 251, Theorem D shows that the hypothesis of Theorem E are satisfied 
if the set $E$ is a subset of $T$ of Lebesgue measure zero and of type $F_\sigma$ and $G_\delta$.


\end{section}

\begin{section}{Proof of Theorem 1.}

{\it Necessity.} Assume that there exists a Blaschke product $B$ such that it has unrestricted
limits at each point of $T \setminus E$ and does not have radial limits at the points of $E$.
First of all, Theorem A implies that $m(E)=0$.
To prove that $E$ is closed, 
let $\{a_n\}  \subset \Delta$ be the set (sequence) of zeros of $B$.
We complete the proof in two steps (propositions).

\underline {Step 1.}  If the unrestricted
limit of $B$ at $z_0 \in T$ exists and is equal to $d$, then $|d|=1$.
 
Indeed, by Theorem B in each neighborhood of $z_0$ on $T$ there exists a point
$w$ at which the radial limit $B(w)$ of $B$ exists and has modulus $1$.   
This implies that $|d|=1$. 

\underline {Step 2.} The set $T \setminus E$ is open.

Indeed, if $z_0 \in T \setminus E$, then $z_0 \notin \overline {\{a_n\}} $ because otherwise
the unrestricted limit of $B$ at $z_0$ would be equal to $0$, which contradicts to Step 1. 
Thus, there exists an open interval $I \subset T$ such that $z_0 \in I$ and $I$ is disjoint of the (closed) set 
$\overline {\{a_n\}}$. Then by a known theorem (see \cite{Hoff}, p. 68) the Blaschke product
 $B$ is analytic at the points of $I$. Thus $B$ is continuous on $I$ and so it possesses unrestricted limits 
 at the points of $I$. This means that $I \subset T \setminus E$ and thus  $T \setminus E$ is open.
Hence $E$ is closed. The proof of necessity is completed.

{\it Sufficiency.} Let $E$ be a closed set of measure zero on $T$. Since $E$ is closed, it is of type $F_\sigma$
and $G_\delta$. Thus for $E$ the hypothesis of Theorem D are satisfied and (by Theorem D) there exist Blaschke 
products $B_0$ and $B_1$ with properties (i) and (ii) listed in Theorem D.  

This means that the hypothesis of Theorem E indicated in the second sentence of its formulation are satisfied.  
Note that the hypothesis that there exists an analytic function $f_1$ in the unit ball of $H^\infty$, $f_1 \not \equiv 1,$ 
such that $\lim_{r \rightarrow 1}f_1(re^{it})=1$ for $e^{it} \in E$, is satisfied even regardless of Theorem D as 
such $f_1$ exists already in the disc algebra by a theorem of Fatou (cf. \cite{Hoff}, p. 81). 

By the sufficiency part of Theorem C, there exists a function $g \in H^\infty$ which has
no radial limit at each point of $E$ but has unrestricted limit at each point of $T \setminus E$ (obviously,
$g$ can be extended continuously on the open set $T\setminus E$). By dividing to a constant if needed, we may assume that the 
function $g$ is in 
the unit ball of $H^\infty$ and thus one can apply for it Theorem E.
Then, by Theorem E, there exists a Blaschke product $I$ that extends analytically to $T \setminus E$, such that for $e^{it} \in E$ we have
 
$$\lim_{r \rightarrow 1} [I(re^{it}) -  g(re^{it})]=0.$$  

 \vspace{0.25 cm}

Since $\lim_{r \rightarrow 1} g(re^{it})$ does not exist for $e^{it} \in E,$  $\lim_{r \rightarrow 1} I(re^{it})$ does not exist for $e^{it} \in E.$
Since $I$ extends analytically to $T \setminus E$, it has unrestricted limits on $T \setminus E$, and the Blaschke product
$I$ has the needed properties. 
This completes the proof of the sufficiency. 

The theorem is proved.  

\end{section}

\begin{minipage}[t]{6.5cm}
Arthur A. Danielyan\\
Department of Mathematics\\ 
and Statistics\\
University of South Florida\\
Tampa, Florida 33620\\
USA\\
{\small e-mail: adaniely@usf.edu}\\
\end{minipage}

\begin{minipage}[t]{6.5cm}
Spyros Pasias\\
Department of Mathematics\\ 
and Statistics\\
University of South Florida\\
Tampa, Florida 33620\\
USA\\
{\small e-mail: spyrospasias@usf.edu}\\
\end{minipage}

\end{document}